\documentclass[11pt,a4paper]{article}
\usepackage[T1]{fontenc}
\usepackage[utf8]{inputenc}
\usepackage{lmodern,microtype}
\usepackage[margin=1.15in]{geometry}
\usepackage{amsmath,amssymb,amsthm,mathtools,booktabs,array}
\usepackage{xcolor}
\usepackage[colorlinks=true,linkcolor=blue!55!black,citecolor=blue!55!black,urlcolor=blue!55!black]{hyperref}
\hypersetup{pdftitle={Quadratic Rational Transformations and Binomial Defining Equations over Finite Fields},
  pdfkeywords={irreducible binomials, finite fields, rational transformations, quadratic character sums, projective linear groups}}
\allowdisplaybreaks[2]
\newtheorem{theorem}{Theorem}[section]
\newtheorem{proposition}[theorem]{Proposition}
\newtheorem{lemma}[theorem]{Lemma}
\newtheorem{corollary}[theorem]{Corollary}
\theoremstyle{definition}
\newtheorem{definition}[theorem]{Definition}
\newtheorem{example}[theorem]{Example}

\newcommand{\F}{\mathbb F}
\newcommand{\I}{\mathcal I}
\newcommand{\E}{\mathcal E}
\newcommand{\PGL}{\operatorname{PGL}_2}
\newcommand{\Nm}{\operatorname{N}}
\title{Quadratic Rational Transformations and Binomial Defining Equations over Finite Fields}
\author{Hongfeng Wu$^{1}$ and Li Zhu$^{2}$
\setcounter{footnote}{-1}
\footnote{E-Mail addresses: whfmath@gmail.com (H. Wu), lizhumath@pku.edu.cn (L. Zhu).}
\\[0.5ex]
\small $^{1}$College of Science, North China University of Technology, Beijing, China
\\
\small $^{2}$School of Mathematical Sciences, Guizhou Normal University, Guiyang, China}

\date{}

\begin{document}
\maketitle
\begin{abstract}
Let $q=p^\ell\equiv1\pmod4$ and $n=2^h\ge8$. We count monic irreducible polynomials of degree $n$ over $\F_q$ whose roots can be transformed, by a rational function of degree at most two, into roots of irreducible binomials. A rigidity theorem shows that any two such rational functions at a fixed source root differ by a nonzero scalar multiplication or a scaled inversion. This determines the multiplicities in the count. The number of accessible polynomials is $q^3(q^2-1)/4$ plus an explicit sum of squares; it is nondecreasing in $h$ and eventually constant. Equality with the lower bound holds exactly when $p\equiv5\pmod8$ or $2^h\mid p+1$. For $q\equiv5\pmod8$ we also count split quadratic pullbacks of each binomial by branch type and classify the $\PGL(\F_q)$-orbits in the accessible set, with explicit representatives over $\F_5$.
\end{abstract}

\section{Introduction}
\label{sec:intro}
An irreducible binomial $Y^n-a$ provides a defining equation for $\F_{q^n}$ with only two nonzero terms. Given a generator $\alpha$ of this extension, we ask whether a rational function of degree at most two, defined over $\F_q$, takes $\alpha$ to a root of such a binomial. Our aim is to count the monic irreducible polynomials having this property. Different rational functions and different target binomials may produce the same source polynomial, so a count of pullbacks alone does not answer the question.

Quadratic rational pullbacks of irreducible polynomials have been considered in several settings, including the enumerations of Ahmadi~\cite{Ahmadi} and Mattarei--Pizzato~\cite{MattareiPizzato}. Rational transformations and the associated projective actions also occur in~\cite{GowMcGuire,Reis,Schulz}. Here both the quadratic map and the irreducible binomial target vary. The relevant objects are the distinct degree-$n$ factors of their split pullbacks.

We work with $q=p^\ell\equiv1\pmod4$ and $n=2^h\ge8$. The targets are precisely the binomials $Y^n-a$ with $a$ a nonsquare in $\F_q$. Theorem~\ref{thm:rigidity} proves that, for a fixed source root, all binomializing maps of degree at most two are obtained from any one of them by scalar multiplication or scaled inversion. Thus every accessible root has the same number of representations, and the linear and genuinely quadratic strata are disjoint.

Counting split quadratic pullbacks by their branch pairs then yields an exact formula for the number of accessible polynomials. Theorem~\ref{thm:defect} expresses this number as $q^3(q^2-1)/4$ plus a sum of squares and proves monotonicity and stability as $h$ varies. Theorem~\ref{thm:attainment} characterizes equality with the lower bound by the arithmetic condition $p\equiv5\pmod8$ or $2^h\mid p+1$. For $q\equiv5\pmod8$, the character sums also give the numbers of splitting branch pairs of each type for an individual binomial. In this congruence class we determine the $\PGL(\F_q)$-orbit lengths and multiplicities of accessible polynomials, including explicit representatives over $\F_5$.

Section~\ref{sec:setup} fixes the notation. Rigidity and the general counting formulas are established in Sections~\ref{sec:rigidity} and~\ref{sec:characters}; Sections~\ref{sec:closed} and~\ref{sec:orbits} treat branch-type refinements and projective orbits, respectively.

\section{Defining equations and restricted rational accessibility}
\label{sec:setup}
Let $k=\F_q$ and $n\ge2$. Fix an algebraic closure $\overline{k}$ of $k$ and regard $L=\F_{q^n}$ as the unique subfield of $\overline{k}$ with $q^n$ elements. Write
\[
\I_q(n)=\{f\in k[X]: f\text{ is monic and irreducible, }\deg f=n\}.
\]
Any root $\alpha\in L$ of $f\in\I_q(n)$ generates $L$ over $k$. For any $\beta\in L$, write $m_\beta\in k[Y]$ for its monic minimal polynomial over $k$. Thus $m_\beta\in\I_q(n)$ if and only if $k(\beta)=L$.

\begin{definition}\label{def:accessible}
For $r\in\{1,2\}$, let $\E_q(n;2,r)$ consist of the polynomials $f\in\I_q(n)$ for which there are a root $\alpha\in L$ of $f$, a nonconstant rational function $R=A/B\in k(X)$ in lowest terms, and $a\in k^\times$ such that
\[
 \max\{\deg A,\deg B\}\le r,\qquad B(\alpha)\ne0,\qquad
 m_{R(\alpha)}(Y)=Y^n-a.
\]
\end{definition}
The definition is independent of the chosen root $\alpha$, since Frobenius commutes with evaluation of rational functions over $k$.

For $q\equiv1\pmod4$ and $n=2^h$, the binomial irreducibility criterion (see~\cite{HeymanShparlinski}) states that $Y^n-a$ is irreducible if and only if $a$ is a nonsquare in $k^\times$. Write $\chi_E$ for the quadratic character of a finite field $E$ of odd order, extended by $\chi_E(0)=0$, and put $\chi_q=\chi_k$. If $n=2^h\ge8$, the target family is
\begin{equation}\label{eq:targets}
 \mathcal T=\{Y^n-a:a\in k^\times,\ \chi_q(a)=-1\},
 \qquad |\mathcal T|=\frac{q-1}{2}.
\end{equation}

To define the polynomial pullback of $g\in\I_q(n)$ by $R$, take the unique reduced expression $R=A/B$ with $B$ monic and put
\begin{equation}\label{eq:pullback}
 g_R(X)=B(X)^n g(A(X)/B(X))\in k[X].
\end{equation}
For $f\in\I_q(n)$ with a chosen root $\alpha\in L$, one has
\begin{equation}\label{eq:divcriterion}
 f\mid g_R\quad\Longleftrightarrow\quad
 B(\alpha)\ne0\text{ and }g(R(\alpha))=0.
\end{equation}
Indeed, if $B(\alpha)=0$, then $g_R(\alpha)=A(\alpha)^n\ne0$, since $A$ and $B$ are coprime and $g$ is monic. Whenever $g(R(\alpha))=0$, the element $R(\alpha)$ has degree $n$ over $k$ and generates $L$.

If $\deg R=2$, then $\deg g_R=2n$. Indeed, the coefficient of $X^{2n}$ is the homogenization of $g$ evaluated at $([X^2]A,[X^2]B)$. These coordinates represent $R(\infty)\in\mathbf P^1(k)$, and the homogenization has no $k$-rational root. In odd characteristic the extension $k(X)/k(R)$ is separable quadratic, and its nonidentity automorphism is a linear fractional transformation over $k$.

At degree four, complete accessibility can occur: over $\F_5$ every irreducible quartic is accessible to $Y^4-2$ by a rational function of degree at most two. For arbitrary roots $\alpha$ of an irreducible quartic and $\beta$ of $Y^4-2$ in $\F_{5^4}$, the equation $A(\alpha)=\beta B(\alpha)$ with $\deg A,\deg B\le2$ imposes four homogeneous linear conditions on six coefficients. In a nonzero solution $B(\alpha)\ne0$, since otherwise $A(\alpha)=0$ and the degree bounds would imply $A=B=0$. The resulting map $A/B$ is nonconstant, and consequently
\[
 |\E_5(4;2,2)|=|\I_5(4)|=\frac{5^4-5^2}{4}=150.
\]

\section{Rigidity of low-degree binomializations}
\label{sec:rigidity}
Assume for the rest of the paper that
\begin{equation}\label{eq:standing}
 q\equiv1\pmod4,\qquad n=2^h\ge8.
\end{equation}
Let $S\subseteq L^\times$ be the set of roots of all binomials in~\eqref{eq:targets}. Since the distinct polynomials in $\mathcal T$ have disjoint sets of $n$ roots, $|S|=n|\mathcal T|=n(q-1)/2$. Put
\[
 \Gamma=\{Y\mapsto cY,\ Y\mapsto c/Y:c\in k^\times\},
 \qquad |\Gamma|=2(q-1).
\]
Multiplication by $c\in k^\times$ changes the parameter $a$ to $c^na$, and inversion changes it to $a^{-1}$. Both preserve the nonsquare class, so $\Gamma$ preserves $S$.

\begin{lemma}\label{lem:power}
Fix $\beta\in S$, with $\beta^n=a\in k^\times$ nonsquare. Every $\beta'\in S$ has the form
\[
 \beta'=c\beta^e,\qquad c\in k^\times,\quad e\in\{1,3,\ldots,n-1\}.
\]
\end{lemma}
\begin{proof}
Both $\beta$ and $\beta'$ belong to the cyclic subgroup
$T=\{z\in L^\times:z^{n(q-1)}=1\}$. The coset $\beta k^\times$ in $L^\times/k^\times$ has order dividing $n$. If it had order $d<n$, then $\beta^d\in k$ and $\beta$ would satisfy $X^d-\beta^d\in k[X]$, contrary to $[k(\beta):k]=n$. The subgroup $\langle\beta,k^\times\rangle$ is the disjoint union of the $n$ cosets $\beta^j k^\times$, $0\le j<n$. It therefore has $n(q-1)$ elements and equals $T$, since $T$ has at most $n(q-1)$ elements. Hence $\beta'=c\beta^e$ for some $e$ modulo $n$. Since $(\beta')^n=c^na^e$ is nonsquare and $c^n$ is a square, $e$ must be odd.
\end{proof}

\begin{theorem}[Low-degree target rigidity]\label{thm:rigidity}
Let $\alpha\in L$ have degree $n$ over $k$. Suppose $R_1,R_2\in k(X)$ are nonconstant rational functions of degree at most two, both defined at $\alpha$, and
\[
 R_1(\alpha),R_2(\alpha)\in S.
\]
Then, as identities in $k(X)$, either
\[
 R_2=cR_1\qquad\text{or}\qquad R_2=c/R_1
\]
for some $c\in k^\times$.
\end{theorem}
\begin{proof}
Put $\beta=R_1(\alpha)$ and $a=\beta^n$, so $a\in k^\times$ is a nonsquare. By Lemma~\ref{lem:power}, write $R_2(\alpha)=c\beta^e$, where $c\in k^\times$ and $e$ is odd with $1\le e<n$. Let $R_i=A_i/B_i$ be reduced expressions. The resultant with respect to $X$ of $A_1(X)-UB_1(X)$ and $A_2(X)-VB_2(X)$ is a nonzero polynomial $H(U,V)$ of degree at most two in each variable, satisfying $H(R_1,R_2)=0$. It is nonzero because, over $k(U,V)$, any root $\xi$ of the first polynomial is algebraic over $k(U)$: if $B_2(\xi)=0$, coprimality gives $A_2(\xi)\ne0$; otherwise $A_2(\xi)/B_2(\xi)$ cannot equal the independent variable $V$.

Write $H(U,V)=\sum_{0\le i,j\le2}h_{ij}U^iV^j$. Substituting $(U,V)=(\beta,c\beta^e)$ and using the basis $1,\beta,\ldots,\beta^{n-1}$ shows that a nonzero $H$ requires a collision
\[
 i+ej\equiv i'+ej'\pmod n,\qquad (i,j)\ne(i',j'),\quad
 0\le i,i',j,j'\le2.
\]
The case $j=j'$ is impossible, since $|i-i'|\le2<n$. If $|j-j'|=1$, parity forces $|i-i'|=1$ and therefore $e\equiv\pm1\pmod n$. If $|j-j'|=2$, then $i-i'$ is even, and it cannot be zero because $n/2\nmid e$. Hence $|i-i'|=2$ and $e\equiv\pm1\pmod{n/2}$. The only possibilities are
\begin{equation}\label{eq:four-exponents}
 e\in\{1,n/2-1,n/2+1,n-1\}.
\end{equation}

Put $m=n/2\ge4$ and list the residues $i+ej$ for $i=0,1,2$ and $j=0,1,2$. If $e=m+1$, the three rows are
\[
 \{0,1,2\},\quad\{m+1,m+2,m+3\},\quad\{2,3,4\}\pmod{2m};
\]
only the residue $2$ occurs twice, from $(2,0)$ and $(0,2)$. If $e=m-1$, the rows are
\[
 \{0,1,2\},\quad\{m-1,m,m+1\},\quad\{2m-2,2m-1,0\}\pmod{2m};
\]
only $0$ occurs twice, from $(0,0)$ and $(2,2)$. This also holds for $m=4$. Thus in the first exceptional case $H$ can contain only the terms $U^2,V^2$, and $H(R_1,R_2)=0$ forces $(R_2/R_1)^2=c^2a$. In the second case only the constant and $U^2V^2$ terms survive, forcing $(R_1R_2)^2=c^2a$. Both identities are impossible: $a$ is a nonsquare in $k$, hence also in $k(X)$.

It follows that $e=1$ or $e=n-1$. Accordingly, $R_2(\alpha)=cR_1(\alpha)$ or $R_2(\alpha)=ca/R_1(\alpha)$. Clearing denominators gives a polynomial of degree at most four vanishing at $\alpha$. Since $[k(\alpha):k]=n\ge8$, that polynomial is identically zero, as required.
\end{proof}

\begin{corollary}[Exact multiplicity]\label{cor:mult}
If a degree-$n$ generator $\alpha$ admits any binomializing rational function of degree at most two, it admits exactly $2(q-1)$ of them, and they form one free $\Gamma$-orbit by postcomposition. Moreover, no such $\alpha$ admits both a linear fractional binomialization and a genuinely quadratic one.
\end{corollary}
\begin{proof}
Theorem~\ref{thm:rigidity} confines all solutions to the orbit of any one solution. Conversely every element of $\Gamma$ preserves $S$. As $R$ is nonconstant, a projective linear map fixing $R$ as a rational function must be the identity, so the orbit is free. Postcomposition by $\Gamma$ preserves degree.
\end{proof}

The degree bound in Theorem~\ref{thm:rigidity} is sharp for every pair $(q,n)$ in~\eqref{eq:standing}. Indeed, if $\beta^n=a$ with $a$ nonsquare, then $(\beta^3)^n=a^3$ is also nonsquare. Thus $X$ and $X^3$ both carry $\beta$ into $S$, and their different degrees place them in distinct $\Gamma$-orbits.

\section{Quadratic pullbacks and character-sum enumeration}
\label{sec:characters}
Write $G=\PGL(k)$, so $|G|=q(q^2-1)$. In characteristic different from two, every separable quadratic rational map $R$ determines a degree-two covering of projective lines with two distinct geometric branch values. These form a Frobenius-stable unordered pair $D$ in the target line.

\begin{lemma}\label{lem:maps-per-branch}
Let $D$ be a Frobenius-stable unordered pair of distinct geometric target points, and let $Q_D(Y)$ be the monic polynomial whose simple roots are the finite points of $D$. Up to a change of source coordinate, the quadratic maps with branch pair $D$ are the two covers
\[
 Z^2=bQ_D(Y),\qquad b\in k^\times/(k^\times)^2.
\]
There are $|G|$ maps for each $D$. The numbers of rational pairs and conjugate quadratic pairs are $q(q+1)/2$ and $q(q-1)/2$, respectively.
\end{lemma}
\begin{proof}
A separable quadratic extension of $k(Y)$ has the form $k(Y)(\sqrt{u})$. Its ramified places are precisely those at which $u$ has odd valuation. Factoring the numerator and denominator in $k[Y]$ shows that a radicand with branch pair $D$ has the form $bQ_D(Y)$ up to a square in $k(Y)$. Its valuation at infinity is odd exactly when $D$ contains infinity. Two constants give the same extension over $k(Y)$ precisely when their ratio is a square in $k^\times$.

Both extensions have a rational parameter over $k$. For $\deg Q_D=1$ the equation exhibits $Y$ as a quadratic polynomial in $Z$. For $\deg Q_D=2$ the projective closure is a smooth conic, which has a $k$-rational point since $k$ is finite. Projection from that point parametrizes the conic. Choices of source parameter form a single $G$-orbit of maps for each extension, with stabilizer the deck group of order two. The two orbits give $2\cdot |G|/2=|G|$ maps. The counts of branch pairs follow from the $q+1$ rational points of the projective line and the $q^2-q$ points defined over $K=\F_{q^2}$ but outside $k$.
\end{proof}

Let $K=\F_{q^2}$. For $a\in k^\times$ nonsquare and $d\ge1$ define
\begin{equation}\label{eq:T12}
 T_1(a,d)=\sum_{t\in k}\chi_q(t^d-a),\qquad
 T_2(a,d)=\sum_{u\in K}\chi_K(u^d-a).
\end{equation}
For $d=n$ every summand equals $\pm1$, since an irreducible polynomial of degree $n\ge8$ has no root in $k$ or in $K$. For a target $g_a(Y)=Y^n-a$, let $B_a^{\rm rat}$ and $B_a^{\rm quad}$ denote the numbers of branch-value pairs with two rational or two conjugate quadratic points, respectively, for which the quadratic pullback of $g_a$ splits into two degree-$n$ irreducible factors. Set $B_a=B_a^{\rm rat}+B_a^{\rm quad}$.

\begin{proposition}[Branch-pair formula]\label{prop:branch}
For each nonsquare $a\in k^\times$,
\begin{align}
 B_a^{\rm rat}&=\frac{q^2-1+(T_1(a,n)+1)^2}{4},\label{eq:Brat}\\
 B_a^{\rm quad}&=\frac{q^2-2q+T_2(a,n)}4,\label{eq:Bquad}\\
 B_a&=\frac{2q^2-2q-1+(T_1(a,n)+1)^2+T_2(a,n)}4.\label{eq:B}
\end{align}
\end{proposition}
\begin{proof}
Fix a quadratic map $R$ with branch pair $D$ and a root $\beta\in L$ of $g_a$. Every point of $D$ is defined over $K$, so $\beta$ is unramified. By Lemma~\ref{lem:maps-per-branch} the fibre splits over $L$ exactly when $bQ_D(\beta)$ is a square in $L$. Since $n$ is even, every element of $k^\times$ is a square in $L$, and both twists have the same splitting behaviour. Neither point above $\beta$ can be infinity, as $R(\infty)\in\mathbf P^1(k)$. In a split fibre the two points generate $L$ over $k$ because their image does. A Frobenius power taking one to the other would fix $\beta$ and hence all of $L$. The two points therefore give distinct degree-$n$ factors, exhausting the pullback of degree $2n$. In a nonsplit fibre either point has degree two over $k(\beta)=L$ and gives an irreducible pullback of degree $2n$.

For a finite rational branch value $t\in k$, assign the sign
\[
 s(t)=\chi_L(\beta-t)=\chi_q\bigl(\Nm_{L/k}(\beta-t)\bigr)=\chi_q(t^n-a),
\]
and put $s(\infty)=1$. A rational branch pair gives a split fibre precisely when its two signs agree. There are $q+1$ signs whose sum is $T_1(a,n)+1$, whence the number of equal-sign unordered pairs is
\[
 \frac{(q+1)^2+(T_1+1)^2-2(q+1)}4
 =\frac{q^2-1+(T_1+1)^2}4.
\]
For a conjugate pair $\{u,u^q\}\subset K\setminus k$, splitting is detected by
\[
 \chi_L((\beta-u)(\beta-u^q))=\chi_K(u^n-a).
\]
Indeed, take the norm of $Q_D(\beta)$ with $Q_D=(Y-u)(Y-u^q)\in k[Y]$. Multiplication over the conjugates of $\beta$ gives $\Nm_{L/k}\bigl(Q_D(\beta)\bigr)=g_a(u)g_a(u^q)$, and hence
\[
\chi_L((\beta-u)(\beta-u^q))
 =\chi_q\bigl((u^n-a)(u^{qn}-a)\bigr)
 =\chi_K(u^n-a).
\]
The value is unchanged by $u\mapsto u^q$. Each $t\in k$ contributes $+1$ to $T_2$ because every nonzero element of $k$ is a square in $K$. Removing these $q$ terms and taking half of the remaining $+1$ values gives $B_a^{\rm quad}=(q^2-2q+T_2)/4$. Adding proves the final expression.
\end{proof}

\begin{theorem}[Exact accessibility counts]\label{thm:general}
Under~\eqref{eq:standing},
\begin{align}
 |\E_q(n;2,1)|&=\frac{q(q^2-1)}4,\label{eq:linearcount}\\
 |\E_q(n;2,2)\setminus\E_q(n;2,1)|
 &=q(q+1)\sum_{\chi_q(a)=-1}B_a,\label{eq:quadcount}
\end{align}
where $B_a$ is given explicitly by~\eqref{eq:B}.
\end{theorem}
\begin{proof}
There are $n(q-1)/2$ target roots. Pairs $(R,\beta)$ with $R\in G$ and $\beta\in S$ contribute $|G|n(q-1)/2$ source-root incidences, since each linear fractional map has a unique preimage in $L$ of every $\beta\in S$. By Corollary~\ref{cor:mult}, every accessible source root is represented $2(q-1)$ times. Dividing also by the $n$ roots per irreducible source polynomial proves~\eqref{eq:linearcount}.

For quadratic maps, Lemma~\ref{lem:maps-per-branch} gives $|G|\sum_aB_a$ pairs $(R,a)$ whose pullback splits. Each supplies $2n$ source roots, while every accessible source root again occurs $2(q-1)$ times by Corollary~\ref{cor:mult}. Therefore the number of genuinely quadratic-accessible source polynomials is
\[
 \frac{2n|G|\sum_aB_a}{2(q-1)n}
 =q(q+1)\sum_aB_a.
\]
\end{proof}

Let $N=\{a\in k^\times:\chi_q(a)=-1\}$ and, for $d$ a power of two with $d\ge4$, put
\[
 F_d(a)=1+T_1(a,d),\qquad H_d=(k^\times)^d .
\]

\begin{theorem}[A sum of squares]\label{thm:defect}
Let $v_2$ denote the $2$-adic valuation. Suppose $q\equiv1\pmod4$ and write $C_h=|\E_q(2^h;2,2)|$ for $h\ge3$. Then
\begin{equation}\label{eq:squares}
 C_h=\frac{q^3(q^2-1)}4+
     \frac{q(q+1)}4\sum_{a\in N}\bigl(F_{2^h}(a)-F_4(a)\bigr)^2.
\end{equation}
In particular, $C_h\ge q^3(q^2-1)/4$, with equality precisely when $F_{2^h}=F_4$ on $N$. Moreover,
\begin{equation}\label{eq:increments}
 C_{h+1}-C_h=\frac{q(q+1)}4\sum_{a\in N}
          \bigl(F_{2^{h+1}}(a)-F_{2^h}(a)\bigr)^2 .
\end{equation}
Thus $C_h$ is nondecreasing in $h$ and is constant for $h\ge\max\{3,v_2(q-1)\}$.
\end{theorem}
\begin{proof}
For $y\in K$ set $M(y)=\sum_{a\in N}\chi_K(y-a)$. If $c\in(k^\times)^2$, then $cN=N$ and $\chi_K(c)=1$, so $M(cy)=M(y)$. Hence on $K^\times$ the function $M$ descends to the cyclic quotient $K^\times/(k^\times)^2$, of order $2(q+1)$. The $2$-part of this order is four. The fourth-power and $2^h$-power maps on the quotient therefore have the same image and the same fibre size. Summing $M(y^{2^h})$ over $y\in K^\times$ and including $y=0$ gives
\begin{equation}\label{eq:aggregate-T2}
 \sum_{a\in N}T_2(a,2^h)=\sum_{a\in N}T_2(a,4).
\end{equation}
The equality in~\eqref{eq:aggregate-T2} holds after summing over $a\in N$; it is not, in general, a termwise identity.

Let $C_a$ be the smooth projective model of $Y^2=X^4-a$. The quartic is squarefree, and its leading coefficient is one, so $C_a$ has genus one and two $k$-rational points at infinity. Its Frobenius trace $t_a=q+1-\#C_a(k)$ is $-F_4(a)$. On the quadratic extension the trace is $t_a^2-2q$. As $\#C_a(K)=q^2+2+T_2(a,4)$, we obtain
\begin{equation}\label{eq:quartic-trace}
 t_a=-F_4(a),\qquad T_2(a,4)=2q-1-F_4(a)^2.
\end{equation}
Substituting \eqref{eq:aggregate-T2} and \eqref{eq:quartic-trace} into \eqref{eq:B}, and using \eqref{eq:linearcount}, \eqref{eq:quadcount} and $|N|=(q-1)/2$, gives
\begin{equation}\label{eq:second-moments}
 C_h=\frac{q^3(q^2-1)}4+
 \frac{q(q+1)}4\sum_{a\in N}\bigl(F_{2^h}(a)^2-F_4(a)^2\bigr).
\end{equation}

For powers of two $4\le d\le D$, the average of $F_D$ on $aH_d$ is $F_d(a)$. Indeed, $\chi_q(-a)=-1$, so
\[
 F_D(a)=\frac{q-1}{|H_D|}\sum_{y\in H_D}\chi_q(y-a).
\]
When $c$ ranges over $H_d$, the element $y/c$ ranges over $H_d$ for every $y\in H_D$. Since $\chi_q(c)=1$, averaging the last display with $a$ replaced by $ca$ gives the claim. In particular $F_d$ is constant on every coset $aH_d$. Summing the identity obtained by expanding $(F_D-F_d)^2$ over these cosets yields
\begin{equation}\label{eq:coset-squares}
 \sum_{a\in N}\bigl(F_D(a)^2-F_d(a)^2\bigr)
   =\sum_{a\in N}\bigl(F_D(a)-F_d(a)\bigr)^2 .
\end{equation}
Apply \eqref{eq:coset-squares} first with $(d,D)=(4,2^h)$ in \eqref{eq:second-moments}, and then with $(d,D)=(2^h,2^{h+1})$ in the difference of two successive instances of that formula. This proves \eqref{eq:squares} and \eqref{eq:increments}. The subgroup $H_{2^h}$ and the fibre size $(q-1)/|H_{2^h}|$ are constant for $h\ge v_2(q-1)$, which gives the asserted stability.
\end{proof}

\begin{theorem}[Attainment of the lower bound]\label{thm:attainment}
Let $q=p^\ell\equiv1\pmod4$, where $p$ is prime, and let $n=2^h\ge8$. Then
\[
 C_h=\frac{q^3(q^2-1)}4
 \quad\Longleftrightarrow\quad
 p\equiv5\pmod8\quad\text{or}\quad 2^h\mid p+1 .
\]
\end{theorem}
\begin{proof}
For $a\in N$ and any power of two $d\ge4$, we have
$F_d(a)=\sum_{x\in k^\times}\chi_q(x^d-a)$, since $\chi_q(-a)=-1$.
Frobenius permutes the summands, and replacing $x$ by $x^{-1}$ gives
\begin{equation}\label{eq:F-symmetries}
 F_d(a^p)=F_d(a),\qquad F_d(a^{-1})=-F_d(a).
\end{equation}
Indeed, $\chi_q(x^{-d}-a)=\chi_q(-a)\chi_q(x^d-a^{-1})$ for $x\ne0$.
Recall from the proof of Theorem~\ref{thm:defect} that $F_n$ is constant on each $H_n$-coset in $N$ and has average $F_4(a)$ on $aH_4$.

Suppose first that $p\equiv5\pmod8$. Write $2^s=\gcd(n,q-1)$, so $s\ge2$.
After choosing a generator of $k^\times$, the $H_n$-cosets in $N$ are indexed by odd exponents modulo $2^s$.
Frobenius multiplies these exponents by $p$. The subgroup generated by $p$ modulo $2^s$ consists of all units congruent to $1$ modulo $4$.
Indeed, $v_2(p-1)=2$, and successive squaring gives $v_2(p^{2^j}-1)=j+2$ for $j\ge0$.
Thus the $H_n$-cosets within each nonsquare $H_4$-coset form one Frobenius orbit. By \eqref{eq:F-symmetries}, $F_n$ is constant there and equals its average $F_4$.
If instead $n\mid p+1$, then $a^{p+1}=(a^{(p+1)/n})^n\in H_n$, so $a^p$ and $a^{-1}$ lie in the same $H_n$-coset.
Equation~\eqref{eq:F-symmetries} gives $F_n(a)=-F_n(a)$ in $\mathbb Z$. Hence $F_n=0$ on $N$, and its average $F_4$ also vanishes.
In both cases \eqref{eq:squares} gives equality.

For the converse, put $M=(q-1)/2=|N|$. Suppose that $D=2^j$ divides $q-1$ and that $8\le D\le n$. Set $m=(q-1)/D$.
We claim that
\begin{equation}\label{eq:digit-obstruction}
 \binom{M}{m}\not\equiv0\pmod p
 \quad\Longrightarrow\quad C_h>\frac{q^3(q^2-1)}4 .
\end{equation}
For $d\in\{4,D\}$, reduce the integer $F_d(a)$ in $k$ and use $\chi_q(z)=z^M$ in $k$. Expanding and summing over $x\in k^\times$ gives $F_d(a)=P_d(a)$ in $k$, where
\[
 P_d(A)=-\!\!\sum_{\substack{0\le i\le M\\(q-1)\mid di}}
              \binom{M}{i}(-A)^{M-i}\ \in k[A].
\]
The terms with $i=0$ cancel in $P_D-P_4$. The smallest positive index for $P_D$ is $m$, whereas for $P_4$ it is $(q-1)/4>m$.
If $\binom{M}{m}$ is nonzero modulo $p$, then $P_D-P_4$ has degree $M-m<M$ and cannot vanish at all $M$ elements of $N$. Thus $F_D(a)\ne F_4(a)$ as integers for at least one nonsquare $a$. Formula~\eqref{eq:squares} gives strict inequality at $j$, and \eqref{eq:increments} carries it to $h$.

It remains to find such a $D$ when neither equality condition holds.
Write $b=(p-1)/2$. The base-$p$ digits of $M$ are all $b$, and in $\F_p[Z]$
\[
 (1+Z)^M=\prod_{i=0}^{\ell-1}(1+Z^{p^i})^b.
\]
Consequently $\binom{M}{m}$ is nonzero modulo $p$ whenever every base-$p$ digit of $m$ is at most $b$.
If $p\equiv1\pmod8$, take $D=8$. Every digit of $m$ equals $(p-1)/8$, so \eqref{eq:digit-obstruction} applies.

The remaining case is $p\equiv3\pmod4$ with $n\nmid p+1$. Here $\ell$ is even. Put $t=v_2(p+1)$ and $D=2^{t+1}$. Then $8\le D\le n$ and $D\mid p^2-1\mid q-1$.
Write $p+1=(2u+1)D/2$. Since
\[
 \frac{p^2-1}{D}=up+(b-u),\qquad 0\le u\le b,
\]
both digits of this integer are at most $b$. Multiplication by $1+p^2+\cdots+p^{\ell-2}$ repeats this pair of digits and gives $m=(q-1)/D$. Equation~\eqref{eq:digit-obstruction} again yields strict inequality.
\end{proof}

\begin{example}\label{ex:seventeen}
The fields $\F_{17}$ and $\F_{25}$ both have order $1$ modulo $8$, but their accessibility counts differ. Over $\F_{17}$,
$H_4=\{1,4,13,16\}$, $H_8=\{1,16\}$, and $H_{16}=\{1\}$. In particular, $F_4(3)=8$ and $F_8(3)=16$. The sums of $F_d(a)^2$ over the eight nonsquares are $512$, $1024$, and $2048$ for $d=4,8,16$, respectively. Hence $C_3=392904$ and $C_h=471240$ for $h\ge4$. Over $\F_{25}$, the fifth-power Frobenius interchanges the two $H_8$-cosets in each nonsquare $H_4$-coset; thus $F_n=F_4$ for all $h\ge3$. Over $\F_{49}$ equality holds only for $h=3$, by Theorem~\ref{thm:attainment}.
\end{example}

\section{Branch-type counts for \texorpdfstring{$q\equiv5\pmod8$}{q congruent to 5 modulo 8}}
\label{sec:closed}
If $q\equiv5\pmod8$, then $v_2(q-1)=2$ and $v_2(q^2-1)=3$, so the power-map images and multiplicities in~\eqref{eq:T12} give, for every $a\in N$,
\begin{equation}\label{eq:4-8}
 T_1(a,n)=\sum_{x\in k}\chi_q(x^4-a),\qquad
 T_2(a,n)=\sum_{x\in K}\chi_K(x^8-a).
\end{equation}
\begin{corollary}[Sizes of the accessibility strata]\label{cor:closed}
Let $q\equiv5\pmod8$ and $n=2^h\ge8$. Then
\begin{align}
 |\E_q(n;2,1)|&=\frac{q(q^2-1)}4,\label{eq:result1}\\
 |\E_q(n;2,2)\setminus\E_q(n;2,1)|&=\frac{q(q^2-1)^2}4,\label{eq:result2}\\
 |\E_q(n;2,2)|&=\frac{q^3(q^2-1)}4.\label{eq:result3}
\end{align}
\end{corollary}
\begin{proof}
By \eqref{eq:4-8}, $F_n=F_4$ on $N$, so \eqref{eq:squares} gives \eqref{eq:result3}. The linear stratum has size \eqref{eq:linearcount}, and subtraction gives \eqref{eq:result2}.
\end{proof}

When $q\equiv5\pmod8$, the branch counts for a fixed target depend on the octic sum in~\eqref{eq:4-8}. The following cancellation identity holds for every $q\equiv1\pmod4$.

\begin{lemma}[Quartic--octic cancellation]\label{lem:cancel}
For every $q\equiv1\pmod4$ and $a\in N$,
\begin{equation}\label{eq:oct-to-quart}
 T_2(a,8)=T_2(a,4).
\end{equation}
In particular, if $q\equiv5\pmod8$ and $n=2^h\ge8$, then
\begin{equation}\label{eq:cancel}
 (T_1(a,n)+1)^2+T_2(a,n)=2q-1.
\end{equation}
\end{lemma}
\begin{proof}
Grouping the octic sum by $u=x^2$ gives
\[
 T_2(a,8)-T_2(a,4)=\sum_{u\in K}\chi_K(u)\chi_K(u^4-a)
                   =\sum_{u\in K}\chi_K(u^5-au).
\]
Take $c\in K$ with $c^2=a$. Then $c^q=-c$, and the norm identity gives
$\chi_K(c)=\chi_q(-a)=-1$. For $P(u)=u^5-au$ and $u\ne0$ one has
\[
 P(c/u)=-\frac{ac}{u^6}\,P(u).
\]
The multiplier $-ac/u^6$ is a nonsquare in $K$: $-a\in k^\times$ is a square in $K$, whereas $c$ is not. The involution $u\mapsto c/u$ has no fixed point and pairs terms of opposite sign on $K^\times$, proving~\eqref{eq:oct-to-quart}. Equations \eqref{eq:4-8} and \eqref{eq:quartic-trace} now give~\eqref{eq:cancel}.
\end{proof}

The parameter $a=3$ over $\F_{17}$ also shows that~\eqref{eq:aggregate-T2} need not hold for individual $a$. Equation~\eqref{eq:quartic-trace} and Lemma~\ref{lem:cancel} give
$T_2(3,4)=T_2(3,8)=-31$, whereas $T_2(3,16)=33$. To verify the latter, the nonzero sixteenth powers in $K=\F_{17^2}$ form the norm-one subgroup $\mu_{18}$, with $16$ preimages each. For $w\in\mu_{18}$, put $s=w+w^{-1}$; then $\Nm_{K/k}(w-3)=10-3s$. The traces $s=\pm2$ occur once, and $s=\pm1,\pm3,\pm4,\pm7$ twice; the quadratic character is negative at $s=1,7,10,13$. Thus eight of the $18$ elements have negative sign, giving $T_2(3,16)=1+16(10-8)=33$.

\begin{corollary}[Counts by branch type]\label{cor:branch-refinement}
Suppose $q\equiv5\pmod8$ and $n=2^h\ge8$. For $a\in N$ one has
\begin{equation}\label{eq:complementary-branches}
 B_a^{\rm rat}=\frac{q^2-1+t_a^2}{4},\qquad
 B_a^{\rm quad}=\frac{q^2-1-t_a^2}{4},\qquad
 B_a=\frac{q^2-1}{2}.
\end{equation}
Each of the two branch counts is independent of $a$.
\end{corollary}
\begin{proof}
By \eqref{eq:4-8}, \eqref{eq:quartic-trace} and Lemma~\ref{lem:cancel},
$T_1(a,n)+1=-t_a$ and $T_2(a,n)=2q-1-t_a^2$. Substitution into \eqref{eq:Brat} and \eqref{eq:Bquad} gives the formulas. The function $F_4$ is constant on each $H_4$-coset. Inversion interchanges the two nonsquare cosets and changes the sign of $F_4$ by \eqref{eq:F-symmetries}. Thus $t_a^2=F_4(a)^2$ is constant on $N$.
\end{proof}
In particular Hasse's bound gives $0\le B_a^{\rm rat}-B_a^{\rm quad}=t_a^2/2\le 2q$.

\section{Projective orbits for \texorpdfstring{$q\equiv5\pmod8$}{q congruent to 5 modulo 8}}
\label{sec:orbits}
The group $G=\PGL(k)$ acts on $\I_q(n)$ by $P\cdot f=m_{P(\alpha)}$, where $\alpha$ is any root of $f$. This is well defined since $\alpha\notin k$ and Frobenius commutes with $P$. If $R(\alpha)$ is a binomial root, then $(R\circ P^{-1})(P(\alpha))=R(\alpha)$. Thus both accessibility strata are $G$-stable.

\begin{theorem}\label{thm:orbits}
Let $q\equiv5\pmod8$, $n=2^h\ge8$, and $G=\PGL(\F_q)$.
The linear fractional stratum is a single orbit of length $|G|/4$.
The genuinely quadratic stratum has $q-1$ orbits of length $|G|/2$
and $(q-1)^2/4$ orbits of length $|G|$. Their stabilizers have orders
two and one, respectively.
\end{theorem}
\begin{proof}
Fix a nonsquare $a\in k^\times$ and put $g_a(Y)=Y^n-a$.
Since $\gcd(n,q-1)=4$, the quotient $k^\times/(k^\times)^n$ is cyclic
of order four. The nonsquare parameters occupy its two classes of order
four, which inversion interchanges. Thus $\Gamma$ acts transitively on the
target polynomials. Its stabilizer of $g_a$ is the group
$H=\{Y\mapsto cY:c^4=1\}$ of order four. Indeed, a reciprocal scaling
would require $c^n=a^2$, whereas the class of $a^2$ has order two in
the quotient. Every source in the linear fractional stratum is consequently
in the $G$-orbit of $g_a$. By Theorem~\ref{thm:rigidity}, a member of $G$
preserving the root set of $g_a$ belongs to $\Gamma$, so its stabilizer
is precisely $H$.

For the genuinely quadratic stratum, it now suffices to consider pullbacks
of $g_a$. For each branch pair $D$ over which the pullback splits, let
$Q_D$ be the monic polynomial with simple roots at the finite points of $D$.
By Lemma~\ref{lem:maps-per-branch}, quadratic maps with branch pair $D$,
up to change of source coordinate, correspond to the two covers
$Z^2=bQ_D(Y)$ indexed by $b\in k^\times/(k^\times)^2$.
Both split above the roots of $g_a$ by Proposition~\ref{prop:branch}.
We thus have $2B_a=q^2-1$ pairs $(D,b)$ by
Corollary~\ref{cor:branch-refinement}.

The $G$-orbits of source factors correspond to the $H$-orbits of these
pairs. To see this, first note that the two degree-$n$ factors belonging
to one cover are exchanged by its deck involution. Indeed, the two points
over a fixed root $\beta$ of $g_a$ belong to distinct Frobenius orbits,
as in Proposition~\ref{prop:branch}. Changing the source coordinate
therefore preserves the $G$-orbit of either factor, and postcomposition
by $H$ preserves the inverse image of the root set of $g_a$.
Conversely, suppose that $P\in G$ takes a factor for $R_1$ to a factor
for $R_2$. For a root $\alpha$ of the first factor,
$R_1(\alpha)$ and $(R_2\circ P)(\alpha)$ are roots of $g_a$.
Theorem~\ref{thm:rigidity} gives $R_2\circ P=\gamma\circ R_1$ for
some $\gamma\in\Gamma$. Since $\gamma$ carries one root of $g_a$ to
another and commutes with Frobenius, it belongs to $H$.

Moreover, the order of the stabilizer of a source factor equals the order
of the stabilizer of its pair $(D,b)$ in $H$. Each scalar stabilizing
$(D,b)$ lifts to exactly two changes of source coordinate, differing by the
deck involution. The lifts permute the two degree-$n$ factors, and exactly
one fixes a chosen factor. Conversely, the preceding rigidity argument
associates to each source stabilizer element a unique stabilizing scalar,
since the rational map is nonconstant.

An element of order four in $H$ can preserve only the unordered pair
$\{0,\infty\}$. Indeed, its square would fix both points. This pair gives
no split pullback, since its signs are $\chi_q(-a)=-1$ and $+1$.
The remaining nonidentity scalar is $Y\mapsto-Y$. Besides
$\{0,\infty\}$, the Frobenius-stable pairs fixed by this scalar are exactly
$\{t,-t\}$ with $t^2\in k^\times$, one for each value of $t^2$.
All $q-1$ pairs give split pullbacks. The element $\beta^2$ has degree
$n/2$ over $k$, since $(\beta^2)^{n/2}=a$ and $[k(\beta):k(\beta^2)]\le2$.
Hence $\beta^2-t^2$ belongs to $\F_{q^{n/2}}^\times$ and is a square
in $L$. Both classes $b$ are fixed, as $Q_D(Y)=Y^2-t^2$ is even.
There are therefore $2(q-1)$ pairs $(D,b)$ with stabilizer of order two,
and the remaining $q^2-1-2(q-1)=(q-1)^2$ have trivial stabilizer.
Dividing by their respective $H$-orbit lengths, two and four, gives
the stated numbers.
\end{proof}

\begin{corollary}\label{cor:five-orbits}
For $q=5$ and every $n=2^h\ge8$, the set $\E_5(n;2,2)$ has nine
projective orbits. The linear fractional orbit has length $30$. The eight
quadratic orbits have the representatives in Table~\ref{tab:maps}, four
of length $120$ and four of length $60$. These account for the $30$
linearly accessible sources and all $4\cdot120+4\cdot60=720$ genuinely
quadratic sources, for a total of $750$.
\end{corollary}
\begin{proof}
Choose $a=2$. The signs on the rational target line are positive at
$\infty,1,2,3,4$ and negative at $0$, so there are ten rational branch
pairs giving split pullbacks. For $u\in\F_{25}\setminus\F_5$, set
$w=u^n\in\mu_3$. The norm of $w-2$ is $1$ when $w=1$ and
$(w-2)(w^{-1}-2)=2$ when $w\ne1$. Thus the conjugate pairs giving split
pullbacks come from $u^n=1$. The kernel of the $n$th power map in
$\F_{25}^\times$ has order eight, and its elements outside $\F_5$ are
the four primitive eighth roots. Their squares are $2$ and $3$, giving the
pairs defined by $Y^2-2$ and $Y^2-3$.

Table~\ref{tab:types} groups the resulting $24$ pairs $(D,b)$ under
$H=\F_5^\times$. Multiplication of the target by $2$ sends $D$ to $2D$
and sends the class $b$ to $b\,2^{-\deg Q_D}$, by monic normalization of
$Q_D(Y/2)$. Thus it exchanges the two classes for a linear $Q_D$ and
preserves them for a quadratic $Q_D$. Within each row of
Table~\ref{tab:types} the branch pairs form a single scalar orbit.
The stabilizer of a chosen $D$ fixes both classes, so representatives
of the two orbits are obtained by taking both classes above that $D$.
For each row of Table~\ref{tab:maps}, direct calculation of the
discriminant in $X$ of $A(X)-YB(X)$ for the displayed maps $R=A/B$
gives $4Q_D(Y)$ and $3Q_D(Y)$, respectively. Their constant factors
are in distinct square classes in $\F_5$. This verifies the
representatives, and Theorem~\ref{thm:orbits} gives the orbit lengths.
\end{proof}

\begin{table}[htbp]
\centering\small
\begin{tabular}{@{}lccc@{}}
\toprule
Branch pair & Number of pairs & Classes per pair & Scalar orbits (length) \\
\midrule
$\{\infty,t\}$, $t\in\F_5^\times$ & $4$ & $2$ & $2$ ($4$) \\
Nonzero finite non-opposite pairs & $4$ & $2$ & $2$ ($4$) \\
$\{1,4\}$, $\{2,3\}$ & $2$ & $2$ & $2$ ($2$) \\
Root pairs of $Y^2-2$, $Y^2-3$ & $2$ & $2$ & $2$ ($2$) \\
\bottomrule
\end{tabular}
\caption{Branch pairs giving split pullbacks of $Y^n-2$ over $\F_5$. The last column lists the scalar orbits of $(D,b)$.}
\label{tab:types}
\end{table}

\begin{table}[htbp]
\centering\small
\renewcommand{\arraystretch}{1.6}
\begin{tabular}{@{}lccr@{}}
\toprule
Branch pair & First map & Second map & Source orbit length \\
\midrule
$\{\infty,1\}$ & $X^2+1$ & $2X^2+1$ & $120$ \\
$\{1,2\}$ & $\dfrac{X^2-2}{X^2-1}$ & $\dfrac{X^2-4}{X^2-2}$ & $120$ \\
$\{1,-1\}$ & $\dfrac{X^2+1}{X^2-1}$ & $\dfrac{X^2+2}{X^2-2}$ & $60$ \\
Roots of $Y^2-2$ & $\dfrac{X^2+2}{2X}$ & $\dfrac{3X}{X^2-2}$ & $60$ \\
\bottomrule
\end{tabular}
\caption{The eight genuinely quadratic orbits over $\F_5$. For each map
$R=A/B$, either degree-$n$ factor of $A^n-2B^n$ represents its source orbit.}
\label{tab:maps}
\end{table}

\begin{example}[A quadratic pullback of a degree-eight binomial]
Over $\F_5$, take $g(Y)=Y^8-2$ and $R(X)=X^2+1$.
A direct factorization gives
\begin{align*}
 g(R(X))&=f_+(X)f_-(X),\\
 f_\pm(X)&=X^8-X^6\pm X^5+X^4\mp2X^3+2X^2+2.
\end{align*}
Both factors are irreducible over $\F_5$, and $f_-(X)=f_+(-X)$.
The deck involution $X\mapsto-X$ exchanges the two factors.
\end{example}

\section*{Data availability}
No datasets were generated or analyzed in this study.

\section*{Declaration of competing interest}
The authors declare no competing interests.

\end{document}